\documentclass[12pt,a4paper]{article}

\usepackage{latexsym}
\usepackage{amsfonts}
\usepackage{amsmath}
\usepackage{amssymb,amsfonts, amsmath, mathrsfs}
\newcommand{\diver}{{{\rm{div}\,}}}

\newcommand{\noin}{\noindent}

\newcommand{\OM}{{\Omega}}

\newcommand{\POM}{{\partial \, \Omega}}

\newtheorem{theorem}{\bf Theorem}[section]

\newtheorem{corollary}[theorem]{\bf Corollary}

\newtheorem{remark}[theorem]{\bf Remark}

\begin{document}

\title{Foliations and  Chern-Heinz inequalities}
\author{J. L. M. BARBOSA \and
G. P. BESSA \and   J. F. MONTENEGRO}
\date{\today}
\maketitle
\begin{abstract} We extend the Chern-Heinz inequalities about mean
curvature  and scalar curvature of graphs of $C^{2}$-functions to
leaves  of transversally oriented codimension one $C^{2}$-foliations
of Riemannian manifolds. That extends partially Salavessa's work on
mean curvature of graphs and generalize  results  of
Barbosa-Kenmotsu-Oshikiri \cite{barbosa-kenmotsu-Oshikiri} and
Barbosa-Gomes-Silveira \cite{barbosa-gomes-silveira} about
foliations of $3$-dimensional Riemannian manifolds  by constant mean
curvature surfaces. These Chern-Heinz inequalities for foliations
can be applied to prove  Haymann-Makai-Osserman inequality (lower
bounds of the fundamental tones of bounded open subsets $\Omega
\subset \mathbb{R}^{2}$ in terms of its inradius) for embedded
tubular neighborhoods of simple  curves of $\mathbb{R}^{n}$.

\vspace{.2cm}

\vspace{.5cm}

\end{abstract}
\section{Introduction}

  Heinz \cite{heinz}, proved that if the graph of
 a $C^{2}$-function $z=z(x,y)$ defined on   $ x^{2}+y^{2}<r^{2}$  has
 mean curvature $H(x,y)$ satisfying $\vert H(x,y)\vert \geq c> 0$ then
$r\cdot c\leq 1$ and if  the Gaussian curvature of the graph
satisfies $K(x,y)\geq c>0$ then $r\cdot c^{1/2}\leq 1$. As observed
by Heinz, these inequalities are a qualitative improvement of
inequalities of the form $r\cdot \varrho (c)\leq 1$, for some
constant $\varrho (c)$ depending on $c$,  implicit in Bernstein's
work \cite{bernstein1}, \cite{bernstein2}.
  Chern \cite{chern}, extended Heinz
inequalities to graphs of $C^{2}$-functions $z=z(x)$ defined on
bounded domains $\OM \subset \mathbb{R}^{n}$ with smooth boundaries
$\partial \OM$, showing that if the mean curvature $H(x)$ of the
graphs satisfy $\vert H(x) \vert \geq c>0$ then $ c\leq
\mbox{vol}_{n-1}(\partial \OM)/\mbox{vol}_{n}(\OM)$ and if the
scalar curvatures satisfy $S(x)\geq c>0$ then
 $ c^{1/2}\leq \mbox{vol}_{n-1}(\partial \OM)/\mbox{vol}_{n}(\OM)$.  These inequalities are
 known nowadays  as the
Chern-Heinz inequalities.

  An immediate  corollary of the Chern-Heinz
inequalities is that an entire graph of a $C^{2}$-function $f :
\mathbb{R}^{n}\to \mathbb{R}$ has constant mean curvature $H$
 if and only if
$H=0$ or constant nonnegative scalar curvature $S\geq 0$ if and only
if $S=0$.

 Salavessa in \cite{salavessa} considered  graphs $G(f)\subset M\times N$ of a smooth maps $f:M\to N$
 and   proved that if the
 graph $G(f)$  has
 parallel mean curvature
vector $H_{G(f)}$ then   for every oriented compact domain $\Omega
\subset M$ one has that $\Vert H_{G(f)}\Vert \leq
\mbox{vol}_{n-1}(\partial \OM)/\mbox{vol}_{n}(\OM)$. Recall that the
Cheeger constant of $M$ is defined by ${\textit
h}(M)=\inf_{\Omega}(\mbox{vol}_{n-1}(\partial
\OM)/\mbox{vol}_{n}(\OM))$, where the infimum is taken over all
 relatively compact subsets $\Omega$  of $M$ with smooth
boundary. In particular, if ${\textit h}(M)=0$  then $G(f)$ is
minimal. If ${\rm dim }(N)=1$ this result is valid regardless the
parallelism of its mean curvature vector.

 A graph of a smooth function $f:\Omega \subset
 M\to \mathbb{R}$ generates by vertical translation a transversally oriented codimension one smooth foliation of
 $\Omega\times \mathbb{R}$, where the leaves are all isometric to the graph.
  In this note we show  that the Chern-Heinz
inequalities can be extended to leaves  of  transversally oriented
codimension-one $C^{2}$-foliations of Riemannian manifolds. Before
we state our main result let us recall that the fundamental tone
$\lambda^{\ast}(\Omega )$ of an open subset $\OM \subset M$ of a
smooth Riemannian manifold  is defined by
\begin{equation}\label{eqTomF}
\lambda^{\ast}(\Omega )= \inf\left\{\frac{\int_{\Omega}\vert \nabla
f\vert^{2}}{\int_{\Omega} f^{2}},\,f\in {H_{0}^{1}(\Omega
)\setminus\{0\}}\right\},
\end{equation}
where $H_{0}^{1}(\Omega)$ is the completion  of
$C^{\infty}_{0}(\Omega )$ w.r.t.  the norm $ \Vert\varphi
\Vert_{\Omega}^2=\int_{\Omega}\varphi^{2}+\int_{\Omega} \vert\nabla
\varphi\vert^2. $ When $\OM$ is a bounded  with piecewise smooth
boundary $\POM$ then $\lambda^{\ast}(\OM)$ coincides with the first
  eigenvalue  $\lambda_{1}(\OM)$ of the Laplacian of $\OM$.

\begin{theorem}\label{theoremFoliation}
Let ${\cal F}$ be a transversely orientable codimension-one
$C^{2}$-foliation of a connected open set $\Omega$ of a Riemannian
manifold $M$. Let $\eta$ be a unit vector field on $\OM$ normal to
the leaves of $\cal F$.  Then
\begin{equation}\label{eqFoliation}
2 \sqrt{\lambda^{\ast}(\Omega)}\geq \inf_{F\in {\cal F}}\inf_{x\in
F} \vert H_{F}\vert(x).
\end{equation}
\end{theorem}  This theorem  has a number of interesting consequences, stated
below as corollaries. It imposes strong restrictions for the
existence of foliations by constant mean curvature hypersurfaces on
open sets with zero fundamental tone or with Ricci curvature bounded
below, corollaries \ref{corFoliation}, \ref{corFoliation2}. On the
other hand if an open set can be foliated by constant mean curvature
smooth hypersurfaces then inequality \ref{corFoliation} gives a
lower bound for the fundamental tone, see Theorem (\ref{thmH-O}).

 Barbosa, Kenmotsu and Oshikiri \cite{barbosa-kenmotsu-Oshikiri},
 have studied $C^{3}$-foliations of  complete manifolds $M$ whose
 leaves have constant mean curvature. They have shown that when $M$ is flat,
 noncompact and all leaves of the foliation have the same mean curvature then
 leaves are minimal. The following corollary extends their result
 to Riemannian manifolds with $\lambda^{\ast}(M)=0$.
\begin{corollary}\label{corFoliation}
Let ${\cal F}$ be transversely oriented codimension-one
$C^{2}$-foliation of a Riemannian manifold $M$ with
$\lambda^{\ast}(M)=0$. If the leaves $F\in {\cal F}$ have the same
constant mean curvature then each leaf is minimal.
\end{corollary}

\begin{remark}The class of smooth Riemannian manifolds $M$ with
  $\lambda^{\ast}(M)=0$ is quite large. It contains all closed,
all complete noncompact with  nonnegative Ricci curvature. In fact,
by Cheng's Comparison Theorem \cite{cheng}  all complete Riemannian
manifolds with asymptotically nonnegative Ricci curvature has zero
fundamental tone. Recall that  $M$ has asymptotically nonnegative
Ricci curvature if $Ric_{M}(x)\geq -\psi(dist_{M}(x_{0},x))$, for a
continuous function $\psi:[0,\infty)\to [0, \infty) $ with
$\lim_{t\to \infty}\psi(t)=0$, $x_{0}\in M$.
\end{remark}

Let ${\cal F}$ be a  transversely oriented codimension-one
$C^{2}$-foliation   of  an open subset $\Omega $
 of the $n$-dimensional complete Riemannian manifold $M$ with Ricci curvature
  bounded below $Ric_{M}\geq (n-1)\kappa$ by complete oriented surfaces with the same constant mean
curvature $H$.
 Barbosa-Kenmotsu-Oshikiri,
(\cite{barbosa-kenmotsu-Oshikiri}, Theorem 3.1)  showed that if $M$
is compact and $\kappa=0$ then the leaves of ${\cal F}$  are totally
geodesic and the Ricci curvature of $M$ is zero in the directions
normal to the leaves. They also showed that if
$M=\mathbb{N}^{n}(\kappa)$, the simply connected space form of
constant sectional curvature $\kappa$, if $\kappa\leq 0$ and  $H\geq
(n-1) \sqrt{-\kappa}$ then $H=(n-1)\sqrt{-\kappa}$
(\cite{barbosa-kenmotsu-Oshikiri}, Theorem 3.8).

As a corollary of Theorem (\ref{theoremFoliation}) we extend Theorem
3.8 of Barbosa-Kenmotsu-Oshikiri to $n$-dimensional Riemannian
manifolds Ricci curvature bounded below ${\rm Ric_{M}}\geq (n-1)k$.

\begin{corollary}\label{corFoliation2}
Let ${\cal F}$ be transversely oriented codimension-one
$C^{2}$-foliation of a complete $n$-dimensional Riemannian manifold
$M$ with Ricci curvature ${\rm Ric_{M}}\geq(n-1)k$. Then
\begin{itemize}
\item[i)]$
2\sqrt{\lambda^{\ast}(\mathbb{N}^{n}(k))}\geq \inf_{F\in {\cal
F}}\inf_{x\in F}\vert H_{F}\vert(x).$
\item[]

\item[ii)]
If $\vert H_{F}\vert \geq c>0$ then $k=-a^{2}$ for some positive $a$
satisfying $(n-1)\cdot a\geq c$.
\end{itemize}
\end{corollary}
We observe here that the proof of  Theorem 1.5 of
\cite{bessa-jorge-oliveira} can be straight forward adapted to give
a version Theorem 3.1 of Barbosa-Kenmotsu-Oshikiri
\cite{barbosa-kenmotsu-Oshikiri} for complete Riemannian manifolds
$M$  nonnegative Ricci curvature. We have the following theorem.
\begin{theorem} \label{thmFoliation2}  Let  ${\cal F}$ be transversely
oriented codimension-one $C^{2}$-foliation  of a complete Riemannian
manifold $M$ with  nonnegative Ricci curvature $Ric_{M}\geq 0$.
Suppose that the leaves are complete oriented hypersurfaces with the
same constant mean curvature $H$. Then $H=0$ and each leaf $F\in
{\cal F}$ is stable. If a leaf $F$ is compact then $F$ is totally
geodesic and the Ricci curvature of $M$ is zero in the normal
directions to $F$.

\end{theorem}
\begin{remark}For $n=3$, Schoen \cite{schoen} showed that a
complete stable minimal surface in a $3$-dimensional Riemannian
manifold with nonnegative Ricci curvature is totally geodesic. Thus
for $n=3$, Theorem \ref{thmFoliation2} fully generalizes Theorem 3.1
of \cite{barbosa-kenmotsu-Oshikiri}.
\end{remark}

\noindent \hspace{.5cm}The graphs $G(f)$ of smooth functions
$f:\Omega\to \mathbb{R}$ generates, by vertical translations,
 transversally oriented smooth foliations
of $\Omega \times \mathbb{R}$. Applying Theorem
(\ref{theoremFoliation}) to such foliations we have the following
generalizations of Chern-Heinz results. The  Corollary
(\ref{salavessa1}) were proved before by Isabel Salavessa,
\cite{salavessa}. We state it here because the proof we provide  is
simpler than hers.
\begin{corollary}[Salavessa]\label{salavessa1}Let $f:\OM\subset M\to\mathbb{R}$ be a $C^{2}$-function
defined on a domain $\Omega$. Consider the product metric on
$M\times \mathbb{R}$. Then
\begin{equation}\label{Chern-Heinz-inequality}\sqrt{\inf \vert H\vert}\leq
2\sqrt{\lambda^{\ast}(\Omega)}.\end{equation}In particular,
\begin{itemize} \item[a)] if $\OM=M$,  $\lambda^{\ast}(M)=0$ and $G(f)$ has constant mean
curvature then the graph is a minimal hypersurface. \item[b)]if
$\OM=M$, $M$ closed and $G(f)$ has constant mean curvature $H$ then
$H=0$,
 $f$ is constant and the graph is an slice of the
 product $\Omega \times \mathbb{R}$.\end{itemize}
\end{corollary}
Corollary (\ref{salavessa1}) has the following version for scalar
curvature.
\begin{corollary}\label{thm-scalar}Let $M$ be an $n$-dimensional Riemannian
manifold with non-positive sectional curvature and $f:\OM\subset
M\to\mathbb{R}$ be a $C^{2}$-function defined on a domain $\Omega$.
Consider the product metric on $M\times \mathbb{R}$ and suppose that
the scalar curvature $S$ of the graph
 $G(f)\subset M\times \mathbb{R}$
satisfies $S\geq 0$. Then
\begin{equation}\label{Chern-Heinz-inequality2}\sqrt{\inf S}\leq
2\sqrt{\lambda^{\ast}(\Omega)}.\end{equation} In particular, if
$\OM=M$,  $\lambda^{\ast}(M)=0$ and the graph $G(f)$ has constant
nonnegative scalar curvature $S\geq 0$ then $S =0.$
\end{corollary}

\noindent \hspace{.5cm}Our last result is an extension of  the
Haymann-Makai-Osserman inequality to tubular neighborhood of simple
curves of $\mathbb{R}^{n}$. Recall that the inradius $\rho(\Omega)$
of a connected open set $\Omega\subset M$ of a Riemannian manifold
is defined as $\rho(\Omega)=\sup\{r>0;\, B_{M}(r)\subset \Omega\}, $
where $B_{M}(r)$ is a ball of radius $r$ of $M$. In \cite{makai},
Makai proved that the fundamental tone $\lambda^{\ast}(\Omega)$ of
 a simply connected bounded domain $\Omega\subset
\mathbb{R}^{2}$ with smooth boundary and inradius  $\rho=\rho
(\Omega)$ was bounded below by $\lambda_{1}(\Omega)\geq
1/4\rho^{2}$. Unaware of Makai's result, Haymann \cite{haymann}
proved years later that $\lambda_{1}(\Omega)\geq 1/900\rho^{2}$.
Osserman \cite{osserman} among other things improved Haymann's
estimate back to $\lambda_{1}(\Omega)\geq 1/4\rho^{2}$. The
Haymann-Makai-Osserman inequality is only known to be valid  for
simply connected bounded open sets $\Omega \subset \mathbb{R}^{2}$
with smooth boundary. Theorem (\ref{theoremFoliation}) can be seen
as (a weak form of) Haymann-Makai-Osserman inequality in higher
dimension. For certain foliated open sets of $\mathbb{R}^{n}$ like
embedded tubular neighborhood of simple smooth curves Theorem
(\ref{theoremFoliation}) can be stated as follows.
\begin{theorem}Let $\gamma:I\subset \mathbb{R}\to \mathbb{R}^{n}$ be a simple smooth curve.
 Let $T_{\gamma}(\rho)$ be an embedded tubular
neighborhood of $\gamma$ with variable  radius. Let $\rho>0$ be its
inradius. Then \begin{equation}\lambda^{\ast}(T_{\gamma}(\rho))\geq
\frac{(n-1)^{2}}{4\rho^2}\end{equation}\label{thmH-O}
\end{theorem}

\section{Foliations and Chern-Heinz Inequalities}
\subsection{Proof of Theorem (\ref{theoremFoliation})}

 Let ${\cal F}$ be a transversally oriented codimension one
$C^{2}$-foliation of a connected  open subset $\OM$ of a Riemannian
manifold $M$. This means that there is a unit vector field $\eta$ on
$\OM$ normal to the leaves. We may suppose that the mean curvature
$H_{F}$ of the leaves $F\in{\cal F}$ are such that $\inf_{F\in {\cal
F}}\inf_{x\in F}\vert H_{F}(x)\vert >0$ otherwise there is nothing
to prove. That means that the mean curvatures of the leaves $H_{F}$,
($\stackrel{\rightarrow}{H_{F}}=H_{F}\cdot\eta$),  does not change
sign in $\Omega$.  Thus we may suppose that the mean curvature
$H_{F}$ is positive and bounded from zero, that is, there is a
positive constant $c_F$ such that the mean curvature $H_F$ computed
with respect to the vector field $\eta$ satisfies $H_F \ge c_F$. The
divergence $\diver \,\eta$  at a point $x\in F$ is given by $
\diver\,\eta(x)= H_F(x) $. On the other hand Bessa and Montenegro
\cite{bessa-montenegro} proved that the fundamental tone of $\OM$ is
bounded below by $\lambda^{\ast}(\OM)\geq c(\OM)^{2}/4$ where $
c(\OM)=\sup_{X}(\inf_{\OM} \, \diver\,X/\Vert X\Vert_{\infty})$ with
the supremum taken over all smooth vector fields $X $ on $\OM$ such
that $\inf_{\OM} \diver X> 0$ and $\Vert X\Vert_{\infty}=\sup_{\OM}
\vert X\vert <\infty$.
 The
estimates for the fundamental tone implies that
\begin{equation}
\lambda^{\ast}(\OM)\geq \left[\inf_{F\in {\cal F}}\inf_{x\in
F}H_{F}(x)\right]^{2}/4\;.
\end{equation}
This is equivalent to (\ref{eqFoliation}).
\subsection{Proof of the corollaries}
The Corollary (\ref{corFoliation}) follows immediately from Theorem
(\ref{theoremFoliation}). The Corollary (\ref{corFoliation2}) is a
direct consequence of Theorem (\ref{theoremFoliation}) and the
Cheng's Comparison Theorem. For if $B_{M}(r) $ is a geodesic ball of
radius $r$ in an $n$-dimensional Riemannian manifold $M$ with Ricci
curvature
 ${\rm Ric}_{M}\geq (n-1)k$ bounded  below, Cheng 's Comparison Theorem says that
$\lambda^{\ast}(B_{M}(r))\leq \lambda^{\ast}(B_{\mathbb{M}(k)}(r))$,
where $B_{\mathbb{M}(k)}(r)$ is  a geodesic ball of radius $r$ in
 the simply connected $n$-dimensional space  $\mathbb{M}(k)$
of constant sectional curvature $k$. In particular
$\lambda^{\ast}(M)\leq \lambda^{\ast}(\mathbb{M}(k))$. Thus, by
inequality \ref{eqFoliation} we have that  $$ \inf_{F\in {\cal
F}}\inf_{x\in F}\vert H_{F}\vert(x)\leq
2\sqrt{\lambda^{\ast}(M)}\leq 2\sqrt{\lambda^{\ast}(\mathbb{M}(k))}.
$$ If $\vert H_{F}\vert \geq c>0$ we have that $
\sqrt{\lambda^{\ast}(\mathbb{M}(k))}\geq c/2$. But
$\lambda^{\ast}(\mathbb{M}(k))\neq 0$ if and only if $k=-a^{2}$ in
which case $\lambda^{\ast}(\mathbb{M}(k))=(n-1)^{2}a^{2}/4$ and
$(n-1)\,a\geq c$.
\subsection{Proof of Theorem \ref{thmFoliation2}}We reproduce here
the proof of Theorem 1.5 of \cite{bessa-jorge-oliveira} with proper
modifications that yields the proof of Theorem \ref{thmFoliation2}.
Let ${\cal F}$ be transversely oriented codimension-one
$C^{2}$-foliation  of a complete Riemannian manifold $M$ with
bounded geometry and nonnegative Ricci curvature $Ric_{M}\geq 0$. By
hypothesis the leaves are complete oriented hypersurfaces with the
same constant mean curvature $H$. By inequality \ref{eqFoliation}
each leaf $F$ is minimal. For each point $p$ of $M$ there is  an
oriented leaf $p\in F_{p}\in {\cal F}$. Let $\nu $ be a continuous
unit vector field normal to $F_{p}$. There is a sequence of leaves
$F_{i}$ of the foliation that converges to $F_{p}$ by the side $\nu$
is pointing. Let $L=\triangle + Ric (\nu) + \vert A\vert^{2}$ the
Jacobi operator on $F_{p}$, where $\triangle$ is the Laplacian on
$F_{p}$, $Ric (\nu)$ is the Ricci curvature of $M$ in the direction
$\nu$ and$\vert A\vert$ is the norm of the second fundamental form
of $F_{p}\subset M$. Take a compact subsets $C\subset C'\subset
F_{p}$ and consider a sequence of compacts $C_{i}\subset F_{i}$
converging to $C'$. We may assume that the only solution of $L(v)=0$
in $C'$ with $v=0$ on $\partial C'$ is $v\equiv 0$. Thus there
exists a function $u\in C^{\infty}(C')$ with $L(u)=1$ in $C'$ and
$u=0$ on $\partial C'$. The mean curvature $H(t)$, $\vert t\vert
<\epsilon$ of the immersions $\psi_{t}:C'\to M$ given by
$\psi_{t}(x)=\exp_{x}(t\,u(x)\nu (x))$ has derivative at $t=0$ given
by $2H'(0)=L(u)=1$ on $C'$, see \cite{ros} Thus for small
$0<t<\epsilon$, $H(t)>0$. If $u$ is positive at some interior point
of $C$ then $\psi_{t} (C')$ has a tangency point with some $C_{i}$
and this is not possible by the maximum principle. Thus $u\leq 0$.
If $u(q)=0$ at an interior point of $C'$ we have that $u\equiv 0$
and this is impossible. Thus $u<0$ in the interior of $C'$ and $u=0$
on $\partial C'$. Setting $w=-u$ we have that $w>0$ in the interior
of $C'$ with $L(w)\leq 0$. In particular $w$ restricted to $C$ is a
positive function. Let $u_{1}$ be the first eigenfunction of $C$,
i.e $L(u_{1})+\lambda_{1}(C)u_{1}=0$. Suppose by contradiction that
$\lambda_{1}(C)<0$. Let $h=w-t\,u_{1}>0$ for small $t$. We have that
$L(h)=L(w)-tL(u_{1})\leq t\lambda_{1}(C)u_{1}<0$ on $C$.  Then
$\triangle h<0$ and $h$ has a minimum in the interior of $C$. By the
maximum principle $h$ is constant. Choosing $t$ in such way that
this minimum is zero we have a contradiction. Thus
$\lambda_{1}(C)\geq 0 $ and $C$ is stable. If $F_{p}$ is not compact
we can exhaust $F_{p}$ by stable compact sets and  $F_{p}$ is
stable. If $F_{p}$ is compact, by a result of Fisher-Colbrie
\cite{fisher-colbrie-schoen} there is a positive function
$g:F_{p}\to \mathbb{R}$ solution to $\triangle g -q g =0$ , $-q=Ric
(\nu) + \vert A\vert^{2}$. Integrating over $F_{p}$ we have that
$\smallint_{F_{p}}\triangle g -q g =\smallint_{F_{p}}-q g=0$. This
implies that $q=0$ thus $ Ric (\nu)=0$ and   $\vert A\vert^{2}=0$.
The stability operator in $F_{p}$ is then $L=\triangle$. If $F_{p}$
is not stable then there exists a function $f:F_{p}\to\mathbb{R}$
with $\smallint_{F_{p}}f=0$ satisfying $\triangle f +
\lambda_{1}(F_{p})f=0$, $\lambda_{1}(F_{p})<0$. Let $D_{f}=\{x\in
F_{p},\,f(x)>0\}$ be a nodal set of $f$. Thus $
\lambda_{1}(D_{f})\leq \smallint_{F_{p}}f\triangle
f/\smallint_{F_{p}}<0$. This contradicts the fact that
$\lambda_{1}(C)>0$ proved before.

\subsection{Graphs}
The proof of item i)  Corollary (\ref{salavessa1}) are immediate
once we have inequality \ref{eqFoliation}. Since a graph generates
an oriented  codimension one $C^{2}$-foliation of $\Omega \times
\mathbb{R}$. The item ii) it is also straight forward since  $M$ is
closed  $\lambda^{\ast}(M)=0$. From inequality \ref{eqFoliation} we
have that $H=0$.  Therefore $G(f)$ besides being compact it is
minimal. Observe that the slices $M\times \{t\}$ are totally
geodesic in $M\times \mathbb{R}$. By the maximum principle the
leaves are $M\times \{t\}$. The proof of Corollary \ref{thm-scalar}
is as follows. We may suppose that $c=\inf S>0$ otherwise there is
nothing to prove. Let $p\in G(f)$ and $\{e_{1}, \cdots e_{n}\}$ be
an orthonormal basis for the tangent space $T_{p}G(f)$ of the graph
$G(f)\subset M\times \mathbb{R}$. The Gauss equation for the plane
generated by $e_i, e_j$ is:
\begin{equation}
K(e_{i},e_{j})=\overline{K}(e_{i},e_{j})+\langle {\rm A
}(e_{i},e_{i}),\,  {\rm A} (e_{j},e_{j})\rangle -\vert {\rm A
}(e_{i},e_{j})\vert^{2}
\end{equation} where $K$ represents the Gaussian curvature of
$G(f)$, $\overline{K}$ the sectional curvature of the ambient space
$\OM \times \mathbb{R}$ and ${\rm A}$ the second fundamental form of
the graph. Adding these equations
\begin{equation}
S(p) = \sum_{i,j}\overline{K}(e_{i},e_{j})+ H^{2} - \Vert {\rm A}
\Vert ^{2}
\end{equation}
Since the sectional curvatures  of $M$ are non-positive one has then
$\overline{K}\leq 0$. It follows that $ S(p) \leq H^2(p) $. Now, if
$S\geq c >0$ then $H\geq \sqrt{c}>0$ and by
(\ref{Chern-Heinz-inequality}) we have
(\ref{Chern-Heinz-inequality2}).

\subsection{Haymann-Makai-Osserman inequality}

Let $\gamma:I\to \mathbb{R}^{n}$ be a simple curve and
$T_{\gamma}(\rho)$ be an embedded  tubular neighborhood of $\gamma$
with inradius $\rho$. Recall that the inradius of an open connected
set $\Omega$ is defined as $\sup\{ r>0,\, B(r)\subset\Omega\}$. The
Haymann-Makai-Osserman inequality says that the fundamental tone of
an  open connected, simply connected, subsets of $\mathbb{R}^{2}$
with smooth boundary has a lower bounded depending  only on the
inradius, $\lambda^{\ast}(\Omega)\geq 1/(2\rho)^{2}$,
$\rho=\rho(\Omega)$ the inradius of $\Omega$.  We have a weak
version of this inequality in higher dimension, since we have it
only for embedded tubular neighborhoods. However, these tubes can
knotted and the curve can be closed. The idea of the proof is
trivial application of inequality \ref{eqFoliation}. An embedded
tube $T_{\gamma}(\rho)$ with inradius  can be foliated by spherical
caps of radius $\rho$. Just run with the center of a ball of radius
$\rho$ along $\gamma$. This gives a foliation whose leaves have
constant mean curvature $(n-1)/\rho$. Thus
$\lambda^{\ast}(T_{\gamma}(\rho))\geq (n-1)^{2}/4\rho^{2}$.

\vspace{1cm}

\noin {\bf Acknowledgments:} The research and writing of this work
was partially supported by CNPq-Brazil.

\end{document}